\documentclass[11pt,letterpaper,reqno]{amsart}
\usepackage{amsmath,amsthm,amsfonts,amssymb,mathtools,algpseudocode}
\usepackage{comment,bookmark,bm,color,float}

\hypersetup{pdfstartview={FitH}}

\newtheorem{Theorem}{{\bf Theorem}}[section]

\newtheorem{Algorithm}[Theorem]{{\bf Algorithm}}

\newtheorem{Lemma}[Theorem]{{\bf Lemma}}

\numberwithin{equation}{section}

\newcommand{\calA}{\mathcal{A}}
\newcommand{\calB}{\mathcal{B}}

\newcommand{\calS}{\mathcal{S}}

\newcommand{\R}{\mathbb{R}}
\newcommand{\diag}{\emph{diag}}
\newcommand{\off}{\text{off}}

\newcommand{\Proj}{\text{Proj}}
\newcommand{\tr}{\text{tr}}

\begin{document}

\title[Trace maximization algorithm for the approximate tensor diagonalization]{Trace maximization algorithm for the approximate tensor diagonalization}

\author{Erna Begovi\'{c}~Kova\v{c}}\thanks{\textsc{Erna Begovi\'{c}~Kova\v{c},
Faculty of Chemical Engineering and Technology, University of Zagreb, Maruli\'{c}ev trg 19, 10000 Zagreb, Croatia}.
\texttt{ebegovic@fkit.hr}}
\author{Ana Perkovi\'{c} (Bok\v{s}i\'{c})}\thanks{\textsc{Ana Perkovi\'{c},
Faculty of Chemical Engineering and Technology, University of Zagreb, Maruli\'{c}ev trg 19, 10000 Zagreb, Croatia}.
\texttt{aboksic@fkit.hr}}

\thanks{This work has been supported in part by Croatian Science Foundation under the project UIP-2019-04-5200.}
\date{\today}

\renewcommand{\subjclassname}{\textup{2020} Mathematics Subject Classification}

\maketitle

\begin{abstract}
In this paper we develop a Jacobi-type algorithm for the approximate diagonalization of tensors of order $d\geq3$ via tensor trace maximization.
For a general tensor this is an alternating least squares algorithm and the rotation matrices are chosen in each mode one-by-one to maximize the tensor trace.
On the other hand, for symmetric tensors we discuss a structure-preserving variant of this algorithm where in each iteration the same rotation is applied in all modes.
We show that both versions of the algorithm converge to the stationary points of the corresponding objective functions.
\end{abstract}

\section{Introduction}

Tensor diagonalization has applications in independent component analysis~\cite{DeLDeMV2001} and signal processing problems like blind source separation, image denoising, etc.~\cite{TPC17,PTC17,Cich15,ComonSP}. The problem has been studied as the orthogonal~\cite{Be21,LUC20,LUC19,LUC18,MMVL08} and non-orthogonal~\cite{TPC17} tensor diagonalization, for structured and unstructured tensors.
In this paper we are interested in the orthogonal tensor diagonalization of a tensor $\calA\in\R^{n\times n\times\cdots\times n}$ of order $d\geq3$. We are looking for the decomposition of the form
\begin{equation}\label{decomp}
\calA=\calS\times_1U_1\times_2U_2\cdots\times_dU_d,
\end{equation}
where $U_m$ are $n_m\times n_m$ orthogonal matrices, $m=1,2,\ldots,d$.

If tensor $\calA$ allows orthogonal diagonalization, then core tensor $\calS$ from the decomposition~\eqref{decomp} is a diagonal tensor. In most of the cases it is not possible to completely diagonlize a tensor using orthogonal transformations. Therefore, a diagonal tensor $\calS$ is not achievable, but we want to maximize its diagonal in a certain way.
In papers~\cite{Be21,LUC20,LUC19,LUC18} authors develop the Jacobi-type algorithms for maximizing the Frobenius norm of the diagonal of $\calS$, that is, maximizing the sum of the squares of the diagonal elements of $\calS$. Here, inspired by the algorithm of Moravitz Martin and Van Loan~\cite{MMVL08}, we design a Jacobi-type algorithm that maximizes the trace of $\calS$.
For a given tensor $\calA\in\R^{n\times n\times\cdots\times n}$ we are looking for its decomposition of the form~\eqref{decomp} such that
$$\tr(\calS)=\sum_{i=1}^d\calS_{ii\cdots i}\rightarrow\max.$$

Using the properties of the mode-$m$ product one can express the core tensor $\calS$ as
$$\calS=\calA\times_1U_1^T\times_2U_2^T\cdots\times_dU_d^T.$$
Then our problem can be written in the following way. Given a tensor $\calA$ we need to find $d$ orthogonal matrices $U_1,U_2,\ldots,U_d$ that maximize the objective function
\begin{equation}\label{problem}
f(U_1,U_2,\ldots,U_d)=\tr(\calA\times_1U_1^T\times_2U_2^T\cdots\times_dU_d^T).
\end{equation}
To solve this problem we develop a Jacobi-type algorithm where each iteration contains $d$ microiterations. In one microiteration we fix $d-1$ matrices and solve the optimization problem for only one matrix, i.e.\@ we optimize in only one mode at a time. This approach is known as alternating least squares (ALS).

Compared to the algorithms where the squares of the diagonal elements are maximized, maximization of the trace leads to a simpler algorithm. Although the trace maximization is not equivalent to the maximization of the Frobenius norm of the diagonal, our numerical examples show that the off-norm of a tensor is decreasing, i.e.\@ the Frobenius norm of the diagonal is increasing, when the trace is increasing. Moreover, we numerically compare the results obtained by our new algorithm with an algorithm that maximizes the squares of the diagonal elements and see that their behaviour is similar.

Apart from the paper~\cite{MMVL08}, trace maximization was addressed in~\cite{ComSor07}. However, none of those papers offers the convergence proof, while here we prove the convergence of our algorithm.
The convergence results are analogous to those from~\cite{Be21,LUC18}.
Since we are maximizing the trace, our objective function is different than the ones in~\cite{Be21,LUC18}, and it is a function of $d$ variables for tensors of order $d$. In particular, we are going to prove that
every accumulation point $(U_1,U_2,\ldots,U_d)$ obtained by our algorithm is a stationary point of the function $f$ defined by~\eqref{problem}.

Moreover, we adapt our trace maximization algorithm to obtain a structure-preserving algorithm for symmetric tensors. Such algorithm will no longer be an ALS algorithm, since we need to optimize over all modes at once, but the convergence theory will be alongside the non-structured ALS algorithm.

In Section~\ref{sec:notation} we shortly introduce the notation. We give a detailed description of the Jacobi-type algorithm for the trace maximization in Section~\ref{sec:algorithm}. In Section~\ref{sec:sym} we study the trace maximization algorithm for symmetric tensors. The convergence proofs are in Section~\ref{sec:cvg}. We end the paper with the numerical results in Section~\ref{sec:numerical}.

\section{Preliminaries and notation}\label{sec:notation}

In this paper we study order-$d$ tensors, $d\geq3$. We use the notation from~\cite{KB09}. Matrices are denoted by capital letters $(A,B,\ldots)$, and tensors by calligraphic letters $(\calA,\calB,\ldots)$. Tensor analogues of matrix rows and columns are called \emph{fibers}. Mode-$m$ fibers are obtained from a tensor by fixing all indices except the $m$th one.
It is often useful to have a matrix representation of a tensor. Mode-$m$ \emph{matricization} of a tensor $\calA\in\R^{n_1\times n_2\cdots \times n_d}$ is a matrix $A_{(m)}\in\R^{n_m\times(n_1\cdots n_{m-1}n_{m+1}\cdots n_d)}$ such that the columns of $A_{(m)}$ are mode-$m$ fibers of $\calA$.

The left and the right multiplication of a matrix $A$ by a matrix $X$ can be generalized as the multiplication of a tensor $\calA$ in mode-$m$. The mode-$m$ \emph{product} of a tensor $\calA\in\R^{n_1\times n_2\times\cdots\times n_d}$ and a matrix $X\in\R^{p\times n_m}$ is a tensor $\calB\in\R^{n_1\times\cdots\times n_{m-1}\times p\times n_{m+1}\times\cdots\times n_d}$,
\begin{equation*}
\calB=\calA\times_m X,
\end{equation*}
such that
\begin{equation*}
B_{(m)}=XA_{(m)}.
\end{equation*}
Two properties that can be derived from the definition of the mode-$m$ product are
\begin{align}
\calA\times_mX\times_nY & =\calA\times_nY\times_mX, \quad m\neq n, \label{product1} \\
\calA\times_mX\times_mY & =\calA\times_m(YX). \label{product2}
\end{align}
The \emph{inner product} of two tensors $\calA,\calB\in\R^{n_1\times n_2\times\cdots\times n_d}$ is given by
$$\langle\calA,\calB\rangle=\sum_{i_1=1}^{n_1}\sum_{i_2=1}^{n_2}\cdots\sum_{i_d=1}^{n_d} a_{i_1i_2\ldots i_d}b_{i_1i_2\ldots i_d},$$
and the \emph{Frobenius norm} of a tensor $\calA$ is defined as
\begin{equation*}
\|\mathcal{A}\|=\sqrt{\langle \mathcal{A},\mathcal{A}\rangle}.
\end{equation*}
This norm is a generalization of the matrix Frobenius norm because
$$\|\mathcal{A}\|=\sqrt{\sum_{i_1=1}^{n_1}\sum_{i_2=1}^{n_2}\cdots\sum_{i_d=1}^{n_d} a_{i_1i_2\ldots i_d}^2}.$$

We say that a tensor $\calA\in\R^{n\times n\times\cdots\times n}$ is \emph{diagonal} if $\mathcal{A}_{i_1 i_2 \ldots i_d}\neq0$ only if $i_1=i_2=\ldots=i_d$.
On the other side, tensor \emph{off-norm} is defined as the Frobenius norm of its off-diagonal part, that is,
$$\off^2(\calA)=\|\calA\|^2-\|\diag(\calA)\|^2.$$

\section{Tensor-trace maximization Jacobi-type algorithm}\label{sec:algorithm}

The idea is to maximize the tensor trace instead of the sum of its squared diagonal entries, which results in a much simpler algorithm. This idea is not new, it has been proposed in the past, but the convergence of such algorithm has only been studied in the (more complicated) symmetric case. Our algorithm described in this section does not impose symmetry. Hence, it may be more attractive because it is much simpler, end therefore, less computationally demanding. However, since it converges to a stationary point of the objective function, it does not always converge to the right solution.

Let $\calA\in\R^{n\times n\times\cdots\times n}$ be an order $d\geq3$ tensor with dimension $n$. We want to find its orthogonal decomposition~\eqref{decomp}
\begin{align*}
\calA & =\calS\times_1U_1\times_2U_2\cdots\times_dU_d, \\
\calS & =\calA\times_1U_1^T\times_2U_2^T\cdots\times_dU_d^T,
\end{align*}
such that the trace of the core tensor $\calS$,
\begin{equation*}
\tr(\calS)=\sum_{i=1}^d\calS_{ii\cdots i},
\end{equation*}
is maximized.
This problem is equivalent to the problem of finding $d$ orthogonal matrices $U_1,U_2,\ldots,U_d$ that maximize the objective function~\eqref{problem}
\begin{equation*}
f(U_1,U_2,\ldots,U_d)=\tr(\calA\times_1U_1^T\times_2U_2^T\cdots\times_dU_d^T)\rightarrow\max.
\end{equation*}

To solve the problem~\eqref{problem} we are using an iterative Jacobi-type algorithm.
In the $k$th iteration we are applying $d$ plane rotations on the underlying tensor $\calA^{(k)}$, one in each mode. We have
\begin{equation}\label{it}
    \calA^{(k+1)}=\calA^{(k)}\times_1R_{1,k}^T\times_2R_{2,k}^T\cdots\times_dR_{d,k}^T, \quad k\geq0, \quad \calA^{(0)}=\calA,
\end{equation}
where $R_{1,k},R_{2,k},\cdots,R_{d,k}\in\R^{n\times n}$ depend on an index pair $(i_k,j_k)$, called pivot pair, and a rotation angle $\phi_k$ as follows,
\begin{equation}\label{rotation}
R_{l,k}=R(i_k,j_k,\phi_k)={\small\left[
    \begin{array}{ccccccccccc}
      1 &  &  &  &  &  &  &  &  &  &  \\
       & \ddots &  &  &  &  &  &  &  &  &  \\
       &  & 1 &  &  &  &  &  &  &  &  \\
       &  &  & \cos\phi_k &  &  &  & -\sin\phi_k &  &  &  \\
       &  &  &  & 1 &  &  &  &  &  &  \\
       &  &  &  &  & \ddots &  &  &  &  &  \\
       &  &  &  &  &  & 1 &  &  &  &  \\
       &  &  & \sin\phi_k &  &  &  & \cos\phi_k &  &  &  \\
       &  &  &  &  &  &  &  & 1 &  &  \\
       &  &  &  &  &  &  &  &  & \ddots &  \\
       &  &  &  &  &  &  &  &  &  & 1 \\
    \end{array}
  \right]
  \begin{array}{l}
     \\
     \\
     \\
     i_k \\
     \\
     \\
     \\
     j_k \\
     \\
     \\
     \\
     \end{array}},
\end{equation}
for $1\leq l\leq d$. The pivot pair $(i_k,j_k)$ is the same for all matrices $R_{l,k}$, $1\leq l\leq d$, but the rotation angle is, in general, different in each mode.

One iteration of our algorithm is made of $d$ microiterations, where $d-1$ variables are held constant, and the remaining one is varied. This approach is called alternating least squares.
The results of $d$ microiterations building the $k$th iteration are denoted by $\calA_l^{(k)}$, $1\leq l\leq d$. We set $\calA_0^{(k)}=\calA^{(k)}$ and $\calA^{(k+1)}=\calA_d^{(k)}$. The microiterations are computed as
\begin{equation}\label{innerit}
\calA_l^{(k)}=\calA_{l-1}^{(k)}\times_1I\cdots\times_{l-1}I\times_lR_{l,k}^T\times_{l+1}I\cdots\times_d I, \quad l=1,\ldots,d.
\end{equation}
Relations~\eqref{innerit} can also be written as matrix products
\begin{equation}\label{matricizationstep}
(\calA_l^{(k)})_{(l)}=R_{l,k}^T(\calA_{l-1}^{(k)})_{(l)}, \quad l=1,\dots,d,
\end{equation}
where each rotation $R_{l,k}$ changes only two rows in the corresponding mode-$l$ matricization $(\calA_{l-1}^{(k)})_{(l)}$.
This scheme is well defined because combining all microiterations~\eqref{innerit} gives the iteration~\eqref{it}. Using the properties of mode-$m$ product~\eqref{product1} and~\eqref{product2} we get
\begin{align*}
\calA^{(k+1)} & = (\calA^{(k)}\times_1R_{1,k}^T\times_2I\cdots\times_dI)\cdots\times_1I\cdots\times_{d-1}I\times_dR_{d,k}^T \\
& = \calA^{(k)}\times_1R_{1,k}^T\times_2R_{2,k}^T\cdots\times_dR_{d,k}^T.
\end{align*}

In the $k$th iteration of the algorithm we have tensor $\calA^{(k)}$ and matrices $U_l^{(k)}$, $1\leq l\leq d$. For the current pivot position $(i_k,j_k)$ we find the rotation matrix $R_{1,k}$. Using $R_{1,k}$ we update the transformation matrix $U_1^{(k)}$ and form the auxiliary tensor $\calA_1^{(k)}$,
$$\calA_1^{(k)}=\calA^{(k)}\times_1R_{1,k}.$$
Since
\begin{align*}
\calA^{(k)}\times_1R_{1,k} & = (\calA\times_1(U_1^{(k)})^T\times_2(U_2^{(k)})^T\cdots\times_d(U_d^{(k)})^T)\times_1R_{1,k}^T \\
& = \calA\times_1(R_{1,k}^TU_1^{(k)})^T\times_2(U_2^{(k)})^T\cdots\times_d(U_d^{(k)})^T,
\end{align*}
it follows that
\begin{equation*}
U_1^{k+1}=U_1^{(k)}R_{1,k}.
\end{equation*}
We repeat the same computation for modes $l=2,\dots d$, one by one, and do the updates
\begin{align*}
\calA_l^{(k)} & =\calA_{l-1}^{(k)}\times_lR_{l,k}^T, \\
U_l^{k+1} & =U_l^{(k)}R_{l,k}.
\end{align*}
We still need to explain how we choose pivot positions $(i_k,j_k)$ and rotations $R_{1,k}$, $k\geq0$.

Regarding the choice of pivot pairs, our algorithm uses cyclic pivot strategies. That means that we go through all possible pivot pairs $(i,j)$, $1\leq i<j\leq n$, in some prescribed order, making a cycle, and then repeat that same cycle until convergence. As we are going to see in Section~\ref{sec:cvg}, the convergence results hold for any cyclic pivot strategy. Still, in order to ensure the convergence we need to set an additional condition --- pivot pair $(i,j)$ must satisfy the condition
\begin{equation}\label{condition}
|\langle\nabla_{U_l}f(U_1,U_2,\ldots,U_d),U_l\dot{R}(i,j,0)\rangle|\geq\eta\|\nabla_{U_l}f(U_1,U_2,\ldots,U_d)\|_2,
\end{equation}
for at least one mode $l$, $1\leq l\leq d$, where $0<\eta\leq\frac{2}{n}$ and $\dot{R}(i,j,0)=\left.\frac{\partial}{\partial\phi}R(i,j,\phi)\right|_{\phi=0}$.
Inequalities of this type are sometimes called \L ojasiewitz gradient inequalities. They are a commonly used tool in proving convergence
of non-linear optimization algorithms, and specifically tensor decomposition algorithms. Similar inequalities were used in e.g.~\cite{Ishteva13} and~\cite{LUC20}.
If a pair $(i,j)$ does not satisfy the condition~\eqref{condition} for any $l$, we move onto the next pair. Even though this condition may seem restrictive, we will show in Section~\ref{sec:cvg} that for every $l=1,\dots,d$ there exists at least one viable pivot pair. Thus, the algorithm will not stop because the condition~\eqref{condition} is not fulfilled.

Now, let us see how the rotation angles are calculated. We fix the index $k$  and assume that the pivot pair is $(i_k,j_k)=(p,q)$, $1\leq p<q\leq n$. We observe an order-$d$ subtensor $\hat{\calA}\in\R^{2\times2\times\cdots\times2}$ of $\calA$. We need to find $2\times2$ rotations $\hat{R}_{l}$, $l=1,\ldots,d$, such that the trace of the subtensor
\begin{equation*}
\hat{\calS}=\hat{\calA}\times_1\hat{R}_{1}^T\times_2\hat{R}_{2}^T\cdots\times_d\hat{R}_{d}^T
\end{equation*}
is maximized. To this end we use mode-$l$ matricizations from~\eqref{matricizationstep}. This gives
\begin{equation}\label{submatricizationstep}
(\hat{\calA_l})_{(l)}=\hat{R}_{l}^T(\hat{\calA}_{l-1})_{(l)}, \quad l=1,\dots,d.
\end{equation}
Since the mode-$l$ matricization is obtained by arranging all mode-$l$ fibers into columns, elements in the same column have all indices the same except for the $l$th one.
Therefore, relation~\eqref{submatricizationstep} can be written as
$$\begin{bmatrix}
    a_{p\ldots p}^{(l)} & \cdots & a_{q\ldots qpq\ldots q}^{(l)} \\
    a_{p\ldots pqp\ldots p}^{(l)} & \cdots & a_{q\ldots q}^{(l)}
    \end{bmatrix}
    =\begin{bmatrix}
     \cos \phi & \sin \phi \\
    -\sin \phi & \cos \phi
    \end{bmatrix}
    \begin{bmatrix}
    a_{p\ldots p}^{(l-1)} & \cdots & a_{q\ldots qpq\ldots q}^{(l-1)} \\
    a_{p\ldots pqp\ldots p}^{(l-1)} & \cdots & a_{q\ldots q}^{(l-1)}
    \end{bmatrix},$$
where in matrices $(\hat{\calA}_{l-1})_{(l)}$ and $(\hat{\calA_l})_{(l)}$ the top right elements have the $l$th index equal to $p$ and the bottom left elements have the $l$th index equal to $q$. In order to maximize the trace of $(\hat{\calA_l})_{(l)}$ we define the function
\begin{align}
g_l(\phi) & =\tr((\hat{\calA_l})_{(l)})=a_{p\ldots p}^{(l)}+a_{q\ldots q}^{(l)} \nonumber \\
& =(\cos\phi a_{p\ldots p}^{(l-1)}+\sin\phi a_{p\ldots pqp\ldots p}^{(l-1)})+(-\sin\phi a_{q\ldots qpq\dots q}^{(l-1)}+\cos\phi a_{q\ldots q}^{(l-1)}). \label{gl}
\end{align}
Setting the derivative of $g_l$ to zero leads to the equation
\begin{align*}
0 ={g_l}'(\phi) & =-\sin\phi a_{p\ldots p}^{(l-1)}+\cos\phi a_{p\ldots pqp\ldots p}^{(l-1)} -\cos\phi a_{q\ldots qpq\ldots q}^{(l-1)}-\sin\phi a_{q\ldots q}^{(l-1)} \\
& = -\sin\phi(a_{p\ldots p}^{(l-1)}+a_{q\ldots q}^{(l-1)})+\cos\phi(a_{p\ldots pqp\ldots p}^{(l-1)}-a_{q\ldots qpq\ldots q}^{(l-1)}).
\end{align*}
By rearranging this equation and dividing it by $\cos\phi$ we get the formula for the tangent of the rotation angle
\begin{equation}\label{tanphi}
\tan\phi=\frac{a_{p\ldots pqp\ldots p}^{(l-1)}-a_{q\ldots qpq\ldots q}^{(l-1)}}{a_{p\ldots p}^{(l-1)}+a_{q\ldots q}^{(l-1)}}.
\end{equation}
This procedure is the same for all $l=1,2,\ldots,d$.

The explicit angles for $\hat{R}_{1}^T,\ldots,\hat{R}_{d}^T$ are not needed in order to perform transformations~\eqref{innerit}. We only need sine and cosine of the corresponding angles. We compute those from \eqref{tanphi} using the transformation formulas
\begin{equation}\label{sincos}
\cos\phi=\frac{1}{\pm \sqrt{1+\tan\phi^2}}, \quad \sin\phi=\frac{\tan\phi}{\pm \sqrt{1+\tan\phi^2}}=\tan\phi\cos\phi.
\end{equation}
We calculate both solutions and take the one that gives a bigger value of the function $g_l$ from~\eqref{gl}.
Notice that the both function values will have the same absolute value but different sign because $\sin\phi_2=-\sin\phi_1$ and $\cos\phi_2=-\cos\phi_1$. Therefore, we can take the angle $\phi_i$, $i=1,2$, that gives a positive value of the function $g_l$.

We sum up this section in Algorithm~\ref{agm:jacobi}.

\begin{Algorithm}\label{agm:jacobi}
\hrule\vspace{1ex}
\emph{Tensor-trace maximization}
\vspace{0.5ex}\hrule
\begin{algorithmic}
\State \textbf{Input:} $\calA^{(0)}\in\R^{n\times n\times\cdots\times n}$, $U_l^{(0)}\in\R^{n\times n}$ orthogonal, $l=1,\dots,d$
\State \textbf{Output:} orthogonal matrices $U_l$, $l=1,\dots,d$
\State $k=0$
\Repeat
\State Choose pivot pair $(i,j)$.
\State $\calA_0^{(k)}=\calA^{(k)}$
\For {l=1:d}
\If {$(i,j)$ satisfies~\eqref{condition} for $l$}
\State Find $\cos\phi_k$ and $\sin\phi_k$ for $R_{l,k}$ using~\eqref{tanphi} and~\eqref{sincos}.
\State $\calA_l^{(k)}=\calA_{l-1}^{(k)}\times_lR_{l,k}^T$
\State $U_l^{k+1}=U_l^{(k)}R_{l,k}$
\EndIf
\EndFor
\State $\calA^{(k+1)}=\calA_d^{(k)}$
\Until{convergence}
\end{algorithmic}
\hrule
\end{Algorithm}

Input arguments in Algorithm~\ref{agm:jacobi} are the initial tensor $\calA^{(0)}$ and the starting app\-ro\-xi\-mations $U_l^{(0)}$, $1\leq l\leq d$.
A simple starting point is to set $\calA^{(0)}=\calA$ and take $U_l^{(0)}$ as identity matrices,
\begin{equation}\label{identinit}
   \calA^{(0)}=\calA, \quad U_l^{(0)}=I,\quad 1\leq l\leq d.
\end{equation}
This identity initialization works very well in most of the cases.

However, identity initialization is, for example, not an option if $\calA$ is an antisymmetric tensor. Recall that a tensor is antisymmetric if
\begin{equation*}
\calA_{\ldots p\ldots q\ldots}=-\calA_{\ldots q\ldots p\ldots},
\end{equation*}
for every pair of indices $(p,q)$.
It is easy to see that this property implies
\begin{equation*}
\calA_{\ldots p\ldots p\ldots}=-\calA_{\ldots p\ldots p\ldots}=0.
\end{equation*}
Hence, only nontrivial elements of an antisymmetric tensor are the ones with all indices different.
Therefore, in equation~\eqref{tanphi} for the tangent of the rotation angle both the numerator and the denominator are equal to zero, and the algorithm fails. That is why we must use a different initialization.
One solution to this problem is to precondition the tensor $\calA$ using the HOSVD (\cite{DeLHOSVD}). We have
\begin{equation*}
\calA=\widetilde{\calS}\times_1\widetilde{U}_1\times_2\widetilde{U}_2\cdots\times_d\widetilde{U}_d,
\end{equation*}
where $\widetilde{U}_l$ are matrices of left singular vectors of matricizations $A_{(l)}$, $1\leq l\leq d$. Then, the HOSVD initialization is given by
\begin{align}
\calA^{(0)} & =\calA\times_1\widetilde{U}_1^T\times_2\widetilde{U}_2^T\cdots\times_d\widetilde{U}_d^T, \nonumber \\
U_l^{(0)} & =\widetilde{U}_l, \quad 1\leq l\leq d. \label{hosvdinit}
\end{align}
In Section~\ref{sec:numerical} we will further discuss these two initializations.

\section{Structure-preserving algorithm for the symmetric tensors}\label{sec:sym}

Algorithm~\ref{agm:jacobi} does not preserve the tensor structure since it applies different rotations in different modes. Still, it can be modified to preserve the symmetry of the starting tensor.

We say that a tensor $\calA$ is symmetric if
\begin{equation*}
\calA_{\ldots p\ldots q\ldots}=\calA_{\ldots q\ldots p\ldots},
\end{equation*}
for any pair of indices $(p,q)$. In other words, entries of a symmetric tensor are invariant to index permutations.
In order to keep the symmetry, transformation matrices $U_1,U_2,\ldots,U_d$ from~\eqref{decomp} should be the same in each mode. That means that now, for a symmetric tensor $\calA$, we are looking for the decomposition of the form
\begin{equation*}
\calA=\calS\times_1U\times_2U\cdots\times_dU,
\end{equation*}
where $U$ is orthogonal.

As we did before, we can write the core tensor $\calS$ as
$$\calS=\calA\times_1U^T\times_2U^T\cdots\times_dU^T.$$
Thus, for a symmetric tensor $\calA$ we need to find the orthogonal matrix $U$ that maximizes the objective function
\begin{equation}\label{symproblem}
f_s(U)=\tr(\calA\times_1U^T\times_2U^T\cdots\times_dU^T).
\end{equation}
Now, in the $k$th iteration of the algorithm we have
\begin{equation}\label{symit}
\calA^{(k+1)}=\calA^{(k)}\times_1R_k^T\times_2R_k^T\cdots\times_dR_k^T, \quad k\geq0, \quad \calA^{(0)}=\calA,
\end{equation}
where $R_k$ is a plane rotation of the form~\eqref{rotation}.
It is interesting to notice that in the matrix case, $d=2$, the trace would remain constant through iterations~\eqref{symit} because
$$\tr(U^TAU)=\tr(A),$$
but that is not the case for the tensors.

Rotations $R_k$ depend on the pivot position and the rotation angle. Pivot positions are chosen in any cyclic order, same as in Algorithm~\ref{agm:jacobi}, with the condition that the pair $(i,j)$ is taken as a pivot pair if it satisfies the condition
\begin{equation}\label{symcondition}
|\langle\nabla f_s(U),U\dot{R}(i,j,0)\rangle|\geq\eta\|\nabla f_s(U)\|_2,
\end{equation}
which is analogue to the condition~\eqref{condition}.

When choosing the rotation angle we now need to consider all modes at once. Hence, this is not an ALS algorithm. Because of that, the formula for the tangent of the rotation angle is more complicated than the from~\eqref{tanphi}. We get a polynomial equation in $\tan\phi$, where the order of the polynomial is equal to the order of the tensor $d$. We derive such equation for $d=3$.

Again, we observe a $2$-dimensional subproblem
$$\hat{\calA}=\hat{\calS}\times_1R\times_2R\times_3R,$$
for a fixed pivot pair $(p,q)$, $1\leq p<q\leq n$, where
$$\hat{\calA}(:,:,1)=\left[
                         \begin{array}{cc}
                           a_{ppp} & a_{pqp} \\
                           a_{qpp} & a_{qqp} \\
                         \end{array}
                       \right], \quad
\hat{\calA}(:,:,2)=\left[
                         \begin{array}{cc}
                           a_{ppq} & a_{pqq} \\
                           a_{qpq} & a_{qqq} \\
                         \end{array}
                       \right],
$$
and
$$R=\left[
      \begin{array}{cc}
        \cos\phi & -\sin\phi \\
        \sin\phi & \cos\phi \\
      \end{array}
    \right].
$$
We choose the angle $\phi$ that maximizes the function
$$g_s(\phi)=\tr(\hat{\calA})=\tr(\hat{\calS}\times_1R^T\times_2R^T\times_3R^T).$$
Using the fact that $\hat{\calA}$ is a symmetric tensor, function $g_s$ can be written as
\begin{align*}
g_s(\phi) & =\cos^3\phi(a_{111}+a_{222})+3\cos^2\phi\sin\phi(a_{112}-a_{122})+ \\
& \qquad +3\cos\phi\sin^2\phi(a_{112}+a_{122})+\sin^3\phi(a_{222}-a_{111}).
\end{align*}
We have
\begin{align*}
0=g_s'(\phi) & =3\cos^3\phi(a_{112}-a_{122})+3\cos^2\phi\sin\phi(2a_{112}+2a_{122}-a_{111}-a_{222})+ \\
& \qquad +3\cos\phi\sin^2\phi(a_{222}-a_{111}-2a_{112}+2a_{122})-3\sin^3\phi(a_{112}+a_{122}).
\end{align*}
Dividing this equation by $-3\cos^3\phi$ we obtain the cubic equation for $t=\tan\phi$,
\begin{align}
& (a_{112}+a_{122})t^3+(a_{111}-a_{222}+2a_{112}-2a_{122})t^2+ \nonumber \\
& \quad +(a_{111}+a_{222}-2a_{112}-2a_{122})t+(a_{122}-a_{112})=0. \label{symtanphi}
\end{align}
We solve this equation for $t$ and calculate $\cos\phi$ and $\sin{\phi}$ using the formulas~\eqref{sincos}. Then we take a real solution of the equation~\eqref{symtanphi} that gives the highest value of the function $g_s$.
This calculation follows the same steps for $d>3$.
From a theoretical point of view, solutions of~\eqref{symtanphi} can be calculated using a rather complicated formula for the roots of the general cubic equation. In practice, we use Matlab function \textit{roots}.

The complete procedure using the identity initialization is given in Algorithm~\ref{agm:symjacobi}.

\begin{Algorithm}\label{agm:symjacobi}
\hrule\vspace{1ex}
\emph{Symmetry-preserving tensor-trace maximization}
\vspace{0.5ex}\hrule
\begin{algorithmic}
\State \textbf{Input:} $\calA\in\R^{n\times n\times\cdots\times n}$ symmetric
\State \textbf{Output:} orthogonal matrix $U$
\State $k=0$
\State $\calA^{(0)}=\calA$
\State $U=I$
\Repeat
\State Choose pivot pair $(i,j)$.
\If {$(i,j)$ satisfies~\eqref{symcondition}}
\State Find $\cos\phi_k$ and $\sin\phi_k$ for $R_k$.
\State $\calA^{(k+1)}=\calA^{(k)}\times_1R_k^T\times_2R_k^T\cdots\times_dR_k^T$
\State $U^{k+1}=U^{(k)}R_k$
\EndIf
\Until{convergence}
\end{algorithmic}
\hrule
\end{Algorithm}

We have observed one intriguing thing. Instead of Algorithm~\ref{agm:symjacobi} for symmetric tensors, one can take its modification where the rotation angle is chosen as the optimal angle in only one (e.g.\@ first) mode, instead of considering all modes at once like it is done in the computation of~\eqref{symtanphi}. The advantage when optimizing the angle in only one mode is that the computation is much simpler, we get a linear equation in $\tan\phi$, the same as in~\eqref{tanphi}. Our convergence proof is valid only if the rotation angle is optimal regarding all modes at once, but the modified algorithm has some interesting properties that can be seen in Figures~\ref{fig:symdiag} and~\ref{fig:symdiag4d} in Section~\ref{sec:numerical}. Note that this modification of Algorithm~\ref{agm:symjacobi} does not lead to the same algorithm as Algorithm~\ref{agm:jacobi}, because it still applies the same rotation in all modes, unlike the Algorithm~\ref{agm:jacobi} when applied to a symmetric tensor.

\section{Convergence proofs}\label{sec:cvg}

In this section we are going to show that Algorithm~\ref{agm:jacobi} and Algorithm~\ref{agm:symjacobi} converge to the stationary points of the objective functions~\eqref{problem} and~\eqref{symproblem}, respectively. The proofs follow the basic idea from the paper~\cite{Ishteva13} that was adopted in~\cite{Be21}.

\subsection{Convergence of Algorithm~\ref{agm:jacobi}}

First, we define the function $\tilde{f}\colon\R^{n\times n}\times\R^{n\times n}\times\cdots\times\R^{n\times n}\to\R$,
\begin{align*}
\tilde{f}(U_1,U_2,\ldots,U_d) &= \tr(\calA\times_1U_1^T\times_2U_2^T\cdots\times_dU_d^T) \\
&= \sum_{r=1}^n\left(\sum_{i_1,\ldots,i_d=1}^n a_{i_1i_2\ldots i_d}u_{i_1r,(1)}u_{i_2r,(2)}\dots u_{i_d r,(d)}\right).
\end{align*}
Function $\tilde{f}$ is the extension of the objective function $f$ from~\eqref{problem} to the set of all square matrices.
We calculate $\nabla_{U_l}\tilde{f}$, $1\leq l\leq d$, element-wise as
\begin{align*}
    \frac{\partial{ \tilde{f}}}{\partial{u_{m r,(l)}} }&=\sum_{i_1,\dots, i_{l-1}, i_{l+1}, \dots, i_d=1}^n a_{i_1\dots i_{l-1} m i_{l+1} \dots i_d} u_{i_1r,(1)}\dots u_{i_{l-1}r,(l-1)}u_{i_{l+1}r,(l+1)}\dots u_{i_d r,(d)}\\
    &= (\calA\times_1U_1^T\cdots\times_{l-1}U_{l-1}^T\times_{l+1}U_{l+1}^T\cdots\times_dU_d^T)_{\underbrace{r\ldots r}_{l-1}m\underbrace{r\ldots r}_{d-l}}.
\end{align*}
Then, $\nabla f$ can be expressed as the projection of $\nabla\tilde{f}$ onto the tangent space at $(U_1,U_2,\ldots,U_d)$ to the manifold $O_n\times O_n\times\cdots\times O_n$, where $O_n$ denotes the set of orthogonal matrices. We have
\begin{align*}
\nabla f(U_1,U_2,\ldots,U_d) &=
\left[
  \begin{array}{ccc}
    \nabla_{U_1}f(U_1,U_2,\ldots,U_d) & \cdots & \nabla_{U_d}f(U_1,U_2,\ldots,U_d) \\
  \end{array}
\right] \\
&=\Proj\left[
  \begin{array}{ccc}
    \nabla_{U_1}\tilde{f}(U_1,U_2,\ldots,U_d) & \cdots & \nabla_{U_d}\tilde{f}(U_1,U_2,\ldots,U_d) \\
  \end{array}
\right] \\
&= \left[
  \begin{array}{ccc}
    U_1\Lambda(U_1) & \cdots & U_d\Lambda(U_d) \\
  \end{array}
\right],
\end{align*}
where
\begin{equation}\label{lambda}
\Lambda(U)\coloneqq\frac{U^T\nabla_U\tilde{f}-(\nabla_U\tilde{f})^TU}{2}.
\end{equation}

Using the operator $\Lambda$ we can simplify the convergence condition~\eqref{condition}.
For $1\leq l\leq d$ we have
$$\|\nabla_{U_l}f(U_1,U_2,\ldots,U_d)\|_2=\|U_l\Lambda(U_l)\|_2=\|\Lambda(U_l)\|_2$$
and
$$\langle\nabla_{U_l}f(U_1,U_2,\ldots,U_d),U_l\dot{R}(i,j,0)\rangle=\langle U_l\Lambda(U_l),U_l\dot{R}(i,j,0)\rangle
=\langle\Lambda(U_l),\dot{R}(i,j,0)\rangle.$$
Since
$$\dot{R}(i,j,0)=
\left[
  \begin{array}{ccccccc}
    0 &  &  &  &  &  &  \\
     & \ddots &  &  &  &  &  \\
     &  & 0 &  & -1 &  &  \\
     &  &  & \ddots &  &  &  \\
     &  & 1 &  & 0 &  &  \\
     &  &  &  &  & \ddots &  \\
     &  &  &  &  &  & 0 \\
  \end{array}
\right]
  \begin{array}{l}
     \\
     \\
     i \\
     \\
     j \\
     \\
     \\
     \end{array},
$$
and $\Lambda(U_l)$ is a skew-symmetric matrix, we get
$$\langle\Lambda(U_l),\dot{R}(i,j,0)\rangle=2\Lambda(U_l)_{ij}.$$
Therefore, the condition~\eqref{condition} can be written as
\begin{equation}\label{conditionlambda}
2|\Lambda(U_l)_{ij}|\geq\eta\|\Lambda(U_l)\|_2.
\end{equation}

Now it is easy to prove that it is always possible to find a pivot pair that satisfies the convergence condition. Lemma~\eqref{lemma1} is a straightforward generalization of Lemma~3.2 from~\cite{Be21}.

\begin{Lemma}\label{lemma1}
For any differentiable function $f\colon O_n\times O_n\times\cdots\times O_n\to\R$, $U_1,U_2,\dots, U_d\in O_n$, and $0<\eta\leq \frac{2}{n}$ it is always possible to find index pairs $(i_{U_l},j_{U_l})$, $1\leq l\leq d$, such that
\begin{equation*}
|\langle\nabla_{U_l}f(U_1,U_2,\ldots,U_d),U_l\dot{R}(i_{U_l},j_{U_l},0)\rangle|\geq\eta\|\nabla_{U_l}f(U_1,U_2,\ldots,U_d)\|_2,
\end{equation*}
where
$\dot{R}(i,j,0)=\left.\frac{\partial}{\partial\phi}R(i,j,\phi)\right|_{\phi=0}$.
\end{Lemma}

\begin{proof}
For any $l=1,\ldots,d$ it is always possible to find an index pair $(p,q)$ such that
\begin{equation*}
|\Lambda(U_l)_{pq}|\geq\frac{1}{n}\|\Lambda(U_l)\|_2.
\end{equation*}
For $\eta=\frac{2}{n}$ we get
\begin{equation*}
2|\Lambda(U_l)_{pq}|\geq\eta\|\Lambda(U_l)\|_2,
\end{equation*}
that is, inequality~\eqref{conditionlambda} is satisfied for $(i_{U_l},j_{U_l})=(p,q)$. As the inequalities~\eqref{conditionlambda} and~\eqref{condition} equivalent, this proves the lemma.
\end{proof}

In Lemma~\ref{lemma2} we show that if $(U_1,U_2,\ldots,U_d)$ is not a stationary point of the function $f$, then applying one step of the Algorithm~\ref{agm:jacobi} to any point in the small enough neighbourhood of $(U_1,U_2,\ldots,U_d)$ would increase the value of $f$. The proof of Lemma~\ref{lemma2} follows the steps of the proof of Lemma~3.4 from~\cite{Be21}.

\begin{Lemma}\label{lemma2}
Let $\{U_l^{(k)}\}_{k\geq 0}$, $1\leq l\leq d$, be the sequences generated by Algorithm~\ref{agm:jacobi}. Let $(\overline{U}_1,\overline{U}_2,\ldots,\overline{U}_d)$ be a $d$-tuple of orthogonal matrices satisfying $\nabla f(\overline{U}_1,\overline{U}_2,\ldots,\overline{U}_d)\neq0$. Then there exist $\epsilon>0$ and $\delta>0$ such that
\begin{equation*}
\|U_l^{(k)}-\overline{U}_l \|_2<\epsilon, \quad \forall l=1,\dots,d,
\end{equation*}
implies
\begin{equation}\label{eq:lemma2}
f(U_1^{(k+1)},U_2^{(k+1)},\ldots,U_d^{(k+1)})- f(U_1^{(k)},U_2^{(k)},\ldots,U_d^{(k)})\geq\delta.
\end{equation}
\end{Lemma}

\begin{proof}
For a fixed iteration $k$ we define $d$ functions $h_k^{(l)}\colon\R\to\R$, $l=1,2,\ldots,d$ as
\begin{align*}
h_k^{(1)}(\phi_1) &= f(U_1^{(k)}R(i_k,j_k,\phi_1),U_d^{(k)},\ldots,U_d^{(k)}), \\
h_k^{(l)}(\phi_l) &= f(U_1^{(k)}R_{1,k},\ldots,U_{l-1}^{(k)}R_{l-1,k},U_l^{(k)}R(i_k,j_k,\phi_l),U_{l+1}^{(k)},\ldots,U_d^{(k)}),\\
h_k^{(d)}(\phi_d) &= f(U_1^{(k)}R_{1,k},U_2^{(k)}R_{2,k},\ldots,U_{d-1}^{(k)}R_{d-1,k},U_d^{(k)}R(i_k,j_k,\phi_d)),
\end{align*}
where $R_{l,k}=R(i_k,j_k,\phi_{U_l,k})$, $2\leq l\leq d-1$.
The rotation angle in Algorithm~\ref{agm:jacobi} is chosen such that
\begin{equation*}
\max_{\phi_l}h_k^{(l)}(\phi_l)=h_k^{(l)}(\phi_{U_l})=f(U_1^{(k)}R_{1,k},\ldots,U_l^{(k)}R_{l,k},U_{l+1}^{(k)},\ldots,U_d^{(k)}), \quad 1\leq l\leq d.
\end{equation*}
Moreover, we know that after each microiteration $l$ in the Algorithm~\ref{agm:jacobi} the value of the objective function $f$ does not decrease, that is,
\begin{align}
f(U_1^{(k+1)},U_2^{(k+1)},\ldots,U_d^{(k+1)}) & \geq f(U_1^{(k+1)},U_2^{(k+1)},\ldots,U_{d-1}^{(k+1)},U_d^{(k)}) \nonumber \\
& \geq\cdots\geq f(U_1^{(k+1)},U_2^{(k)},\ldots,,U_d^{(k)}) \label{lm:inequality} \\
& \geq f(U_1^{(k)},U_2^{(k)},\ldots,U_d^{(k)}). \nonumber
\end{align}
To prove the inequality~\eqref{eq:lemma2} we need at least one sharp inequality in~\eqref{lm:inequality}.

Since $\nabla f(\overline{U}_1,\overline{U}_2,\ldots,\overline{U}_d)\neq0$, we have
\begin{equation*}
\nabla_{U_l}f(\overline{U}_1,\overline{U}_2,\ldots,\overline{U}_d)\neq0,
\end{equation*}
for at least one partial gradient $\nabla_{U_l}f$, $1\leq l\leq d$.
Let us assume that $m$, $1\leq m\leq d$, is the smallest index such that
\begin{equation}\label{lm:gradm}
\nabla_{U_m} f(\overline{U}_1,\overline{U}_2,\ldots,\overline{U}_d)\neq0.
\end{equation}
Then,
\begin{align}
& f(U_1^{(k+1)},U_2^{(k+1)},\ldots,U_d^{(k+1)})-f(U_1^{(k)},U_2^{(k)},\ldots,U_d^{(k)}) \nonumber \\
& \qquad \geq f(U_1^{(k+1)},\ldots,U_m^{(k+1)},U_{m+1}^{(k)},\ldots,U_d^{(k)})-f(U_1^{(k)},U_2^{(k)},\ldots,U_d^{(k)}) \nonumber \\
& \qquad \geq h_k^{(m)}(\phi_m)-h_k^{(m)}(0). \label{lm:fh1}
\end{align}

We need the Taylor expansion of the function $h_k^{(m)}$ around $0$. It is given by
\begin{equation}\label{lm:taylor}
h_k^{(m)}(\phi_m)=h_k^{(m)}(0)+(h_k^{(m)})'(0)\phi_m +\frac{1}{2}(h_k^{(m)})''(\xi)\phi_m^2, \quad 0<\xi<\phi_m.
\end{equation}
Denote $M=\max|(h_k^{(m)})''(\xi)|<\infty$. Then we can write the Taylor expansion~\eqref{lm:taylor} as
\begin{equation}\label{lm:taylorineq}
h_k^{(m)}(\phi_m)-h_k^{(m)}(0) \geq (h_k^{(m)})'(0)\phi_m - \frac{1}{2}M\phi_m^2.
\end{equation}
Therefore, using relations~\eqref{lm:fh1} and~\eqref{lm:taylorineq} we obtain
\begin{equation}\label{lm:fh2}
f(U_1^{(k+1)},U_2^{(k+1)},\ldots,U_d^{(k+1)})-f(U_1^{(k)},U_2^{(k)},\ldots,U_d^{(k)}) \geq (h_k^{(m)})'(0)\phi_m - \frac{1}{2}M\phi_m^2.
\end{equation}

The derivative of $h_k^{(m)}$ is calculated as
\begin{align*}
(h_k^{(m)})'(\phi_m) & =\langle\nabla_{U_m}f(U_1^{(k)}R_{1,k},\ldots,U_{m-1}^{(k)}R_{m-1,k},U_m^{(k)}R(i_k,j_k,\phi_m),U_{m+1}^{(k)},\ldots,U_d^{(k)}), \\
& \qquad \qquad U_m^{(k)}\dot{R}(i_k,j_k,\phi_m)\rangle.
\end{align*}
From the assumption that
\begin{equation*}
\nabla_{U_l} f(\overline{U}_1,\overline{U}_2,\ldots,\overline{U}_d)=0, \quad 1\leq l<m,
\end{equation*}
it follows that the transformations in the first $m-1$ variables did not change the value of $f$,
\begin{align*}
& f(U_1^{(k)}R_{1,k},\ldots,U_{m-1}^{(k)}R_{m-1,k},U_m^{(k)}R(i_k,j_k,\phi_m),U_{m+1}^{(k)},\ldots,U_d^{(k)}) \\
& \quad =f(U_1^{(k)},\ldots,U_{m-1}^{(k)},U_m^{(k)}R(i_k,j_k,\phi_m),U_{m+1}^{(k)},\ldots,U_d^{(k)}).
\end{align*}
Knowing that $R(i_k,j_k,0)=I$, we get the value of $(h_k^{(m)})'$ at $\phi_m=0$,
\begin{equation}\label{lm:der0}
(h_k^{(m)})'(0)=\langle\nabla_{U_m}f(U_1^{(k)},\ldots,U_d^{(k)}),U_m^{(k)}\dot{R}(i_k,j_k,0)\rangle.
\end{equation}
Hence, Lemma~\ref{lemma1} and equation~\eqref{lm:der0} imply
\begin{equation}\label{lm:eta}
|(h_k^{(m)})'(0)|\geq\eta\|\nabla_{U_l}f(U_1^{(k)},U_2^{(k)},\ldots,U_d^{(k)})\|_2.
\end{equation}
It follows from the relation~\eqref{lm:gradm} that there exists $\epsilon>0$ such that
\begin{equation}\label{lm:mu}
\mu\coloneqq\min\{\|\nabla_{U_m}f(U_1,U_2,\ldots,U_d)\|_2 \, \colon\, \|U_m-\overline{U}_m \|_2<\epsilon\}>0.
\end{equation}
Now, relation~\eqref{lm:mu} together with the inequality~\eqref{lm:eta} gives the lower bound on $|(h_k^{(m)})'(0)|$,
\begin{equation}\label{lm:etamu}
|(h_k^{(m)})'(0)|\geq\eta\mu>0.
\end{equation}

Finally, we go back to inequality~\eqref{lm:fh2}. For $\phi_m=\frac{1}{M}(h_k^{(m)})'(0)$, using the relation~\eqref{lm:etamu}, we get
\begin{align*}
& f(U_1^{(k+1)},U_2^{(k+1)},\ldots,U_d^{(k+1)})-f(U_1^{(k)},U_2^{(k)},\ldots,U_d^{(k)}) \\
& \qquad \geq \frac{1}{M}((h_k^{(m)})'(0))^2-\frac{1}{2M}((h_k^{(m)})'(0))^2 \\
& \qquad \geq \frac{1}{2M}\eta^2\mu^2=\delta>0.
\end{align*}
\end{proof}

Using Lemma~\eqref{lemma2} we are going to prove that Algorithm~\ref{agm:jacobi} converges to a stationary point of the objective function.

\begin{Theorem}\label{theorem}
Every accumulation point $(\overline{U}_1,\overline{U}_2,\dots,\overline{U}_d)$ obtained by Algorithm~\ref{agm:jacobi} is a stationary point of the function $f$ defined by~\eqref{problem}.
\end{Theorem}

\begin{proof}
Suppose that $\overline{U}_l$ are the accumulation points of the sequences $\{U^{(j)}_l\}_{j\geq 1}$, $1\leq l\leq d$, generated by Algorithm~\ref{agm:jacobi}. Then there are subsequences $\{U^{(j)}_l\}_{j\in\mathbf{K}_l}$, $\mathbf{K}_l\subseteq\mathbb{N}$, such that
\begin{equation*}
\{U^{(j)}_l\}_{j\in \mathbf{K}_l}\to\overline{U}_l, \quad 1\leq l\leq d.
\end{equation*}

Further on, suppose that
\begin{equation*}
\nabla f(\overline{U}_1,\overline{U}_2,\ldots,\overline{U}_d)\neq 0.
\end{equation*}
Then, for any $\epsilon>0$ there are $K_l\in\mathbf{K}_l$, $1\leq l\leq d$ such that
\begin{equation*}
\|U_l^{(k)}-\overline{U}_l\|_2<\epsilon, \quad \forall l=1,\dots,d,
\end{equation*}
for every $k>K$, $K=\max\{K_l \, \colon\, 1\leq l\leq d\}$. Lemma~\ref{lemma2} implies that
$$f(U_1^{(k+1)},U_2^{(k+1)},\ldots,U_d^{(k+1)})- f(U_1^{(k)},U_2^{(k)},\ldots,U_d^{(k)})\geq\delta,$$
for some $\delta>0$. Therefore, we have
\begin{equation*}
f(U_1^{(k)},U_2^{(k)},\ldots,U_d^{(k)})\to\infty,
\end{equation*}
when $k\to\infty$.

Since $f$ is a continuous function, convergence of $(U^{(j)}_1,U^{(j)}_2,\ldots,U^{(j)}_d)$ implies the convergence of $f(U^{(j)}_1,U^{(j)}_2,\ldots,U^{(j)}_d)$ and we got a contradiction. Hence,
$$\nabla f(\overline{U}_1,\overline{U}_2,\ldots,\overline{U}_d)=0,$$
that is, $(\overline{U}_1,\overline{U}_2,\ldots,\overline{U}_d)$ is a stationary point of the function $f$.
\end{proof}

\subsection{Convergence of Algorithm~\ref{agm:symjacobi}}

To prove the convergence of the Algorithm~\ref{agm:symjacobi} we follow the same scheme as for the Algorithm~\ref{agm:jacobi}. We should keep two things in mind. First, the function that is being maximized by the Algorithm~\ref{agm:symjacobi} is a function of only one variable. Second, unlike the Algorithm~\ref{agm:jacobi}, this is not an ALS algorithm. These two facts will actually simplify the lemmas needed for the proof.

Instead of Lemma~\ref{lemma1} we can now use Lemma~3.1 from~\cite{LUC18}.

\begin{Lemma}[\cite{LUC18}]
For every differentiable function $f_s\colon O_n\to\R$, $U\in O_n$, and $0<\eta\leq \frac{2}{n}$ it is always possible to find index pair $(i,j)$ such that~\eqref{symcondition} holds.
\end{Lemma}

Lemma~\ref{lemma2sym} is similar to Lemma~\ref{lemma2}, but instead of $d$ microiterations we observe one iteration.

\begin{Lemma}\label{lemma2sym}
Let $\{U^{(k)}\}_{k\geq0}$, be the sequence generated by Algorithm~\ref{agm:symjacobi}. For $\overline{U}\in O_n$, let $\nabla f(\overline{U})\neq0$. Then there exist $\epsilon>0$ and $\delta>0$ such that
$\|U^{(k)}-\overline{U}\|_2<\epsilon$
implies
\begin{equation*}
f_s(U^{(k+1)})- f_s(U^{(k)})\geq\delta.
\end{equation*}
\end{Lemma}

\begin{proof}
The proof follows the same reasoning as the proof of Lemma~\ref{lemma2}.

For a fixed iteration $k$ we define the function $h_k\colon\R\to\R$,
\begin{equation*}
h_k(\phi)=f_s(U^{(k)}R(i_k,j_k,\phi)).
\end{equation*}
The rotation angle in Algorithm~\ref{agm:symjacobi} is chosen in such a way that
\begin{equation*}
\max_{\phi}h_k(\phi)=h_k(\phi_k)=f_s(U^{(k)}R(i_k,j_k,\phi_k)=f_s(U^{(k+1)}).
\end{equation*}
Moreover, we have
\begin{equation*}
h_k(0)=f_s(U^{(k)}R(i_k,j_k,0))=f_s(U^{(k)}).
\end{equation*}
Thus,
\begin{equation}\label{lm:symfh1}
f_s(U^{(k+1)})-f_s(U^{(k)})=h_k(\phi_k)-h_k(0).
\end{equation}

We use the Taylor expansion of the function $h_k$ around $0$,
\begin{equation*}
h_k(\phi_k)=h_k(0)+h_k'(0)\phi_k +\frac{1}{2}h_k(\xi)\phi_k^2, \quad 0<\xi<\phi_k.
\end{equation*}
For $M=\max|h_k''(\xi)|<\infty$ it follows from the relation~\eqref{lm:symfh1} that
\begin{equation}\label{lm:symfh2}
f(U^{(k+1)})-f(U^{(k)}) \geq h_k'(0)\phi_k - \frac{1}{2}M\phi_k^2.
\end{equation}
Using Lemma~\ref{lemma2sym} we get
\begin{align}
h_k'(0) & =\langle\nabla f_s(U^{(k)})R(i_k,j_k,0),U^{(k)}\dot{R}(i_k,j_k,0)\rangle \nonumber \\
& =\langle\nabla f_s(U^{(k)}),U^{(k)}\dot{R}(i_k,j_k,0)\rangle \nonumber\\
& \geq\eta\|\nabla f_s(U^{(k)})\| \label{lm:symder0}
\end{align}
Since $\|U^{(k)}-\overline{U}\|_2<\epsilon$, there exists $\epsilon>0$ such that
\begin{equation}\label{lm:symmu}
\mu\coloneqq\min\{\|\nabla f_s(U)\|_2 \, \colon\, \|U-\overline{U} \|_2<\epsilon\}>0,
\end{equation}
and it follows from~\eqref{lm:symder0} and~\eqref{lm:symmu} that
\begin{equation*}
|h_k'(0)|\geq\eta\mu.
\end{equation*}
Then, from the inequality~\eqref{lm:symfh2}, for $\phi_k=\frac{1}{M}(h_k)'(0)$, we get
\begin{equation*}
f(U^{(k+1)})-f(U^{(k)}) \geq\frac{1}{2}\frac{h_k'(0)^2}{M} \geq\frac{1}{2M}\eta^2\mu^2=\delta>0.
\end{equation*}
\end{proof}

Now we can prove the convergence of Algorithm~\ref{agm:symjacobi}.

\begin{Theorem}\label{symtheorem}
Every accumulation point $U$ obtained by Algorithm~\ref{agm:symjacobi} is a stationary point of the function $f_s$ defined by~\eqref{symproblem}.
\end{Theorem}

\begin{proof}
The proof is analogous to the proof of Theorem~\ref{theorem}. Instead of Lemma~\ref{lemma2} it uses Lemma~\ref{lemma2sym}.
\end{proof}

We end this section with the expression for $\nabla f_s$. We define the extension of the objective function $f_s$ as $\tilde{f_s}\colon\R^{n\times n}\to\R$,
\begin{align*}
\tilde{f_s}(U) &=\tr(\calA\times_1U^T\times_2U^T\cdots\times_dU^T)=\sum_{r=1}^n\left(\sum_{i_1,\ldots,i_d=1}^n a_{i_1i_2\ldots i_d}u_{i_1r}u_{i_2r}\cdots u_{i_d r}\right) \\
&= \sum_{r=1}^n\sum_{m=1}^n\sum_{k=1}^d {d \choose k} \left(\sum_{\substack{i_{k+1},\ldots,i_d=1 \\ i_{k+1},\ldots,i_d\neq m}}^n
a_{m\ldots m i_{k+1}\ldots i_d} u_{mr}^k u_{i_{k+1}r}\cdots u_{i_dr}\right).
\end{align*}
Element-wise, the gradient of $\tilde{f_s}$ is given by
\begin{align*}
\frac{\partial{\tilde{f_s}}}{\partial{u_{m r}}}
&= \sum_{k=1}^d {d\choose k} ku_{mr}^{k-1} \left(\sum_{\substack{i_{k+1},\ldots,i_d=1 \\ i_{k+1},\ldots,i_d\neq m}}^n
a_{m\ldots m i_{k+1}\ldots i_d} u_{i_{k+1}r}\cdots u_{i_dr}\right) \\
&= \sum_{k=1}^d {d\choose k} ku_{mr}^{k-1} \left(\calA\times_{k+1}\hat{U}\cdots\times_d\hat{U}\right)_{\underbrace{m\cdots m}_{k}\underbrace{r\cdots r}_{d-k}},
\end{align*}
where $\hat{U}$ is a matrix equal to $U$ in all entries except for the $m$th row where the entries of $\hat{U}$ are equal to zero.
Then, $\nabla f_s$ is the projection of $\nabla\tilde{f_s}$ onto the tangent space at $U$ to the manifold $O_n$. That is,
\begin{equation*}
\nabla f_s(U) =\Proj\nabla\tilde{f_s}(U) =U\Lambda(U),
\end{equation*}
where the operator $\Lambda$ is defined by the relation~\eqref{lambda}.

\section{Numerical experiments}\label{sec:numerical}

In the final section of this paper we present the results of our numerical experiments. All the tests are done in Matlab R2021a.

For both Algorithm~\ref{agm:jacobi} and Algorithm~\ref{agm:symjacobi} we observe two values in each mi\-cro\-i\-te\-ra\-ti\-on --- trace and relative off-norm of a current tensor. The trace is the objective function which is expected to increase in each microiteration and converge to some value. The algorithms stop when the change of the trace after one cycle is less than $10^{-4}$.
The relative off-norm of a tensor $\calA$ is given by
$$\frac{\text{off}(\calA)}{\|\calA\|}.$$
Obviously, relative off-norm of a diagonal tensor is equal to zero.

The algorithms are applied on general random tensors and random tensors that can be diagonalized using orthogonal transformations.
Random tensor entries are drawn from uniform distribution on interval $[0,1]$. Orthogonally diagonalizable tensors are constructed such that we take a diagonal tensor with random uniformly distributed entries from $[0,1]$ on the diagonal and multiply it in each mode with random orthogonal matrices (obtained from QR decomposition of random matrices).

\begin{figure}[h!]
    \centering
    \includegraphics[width=\textwidth]{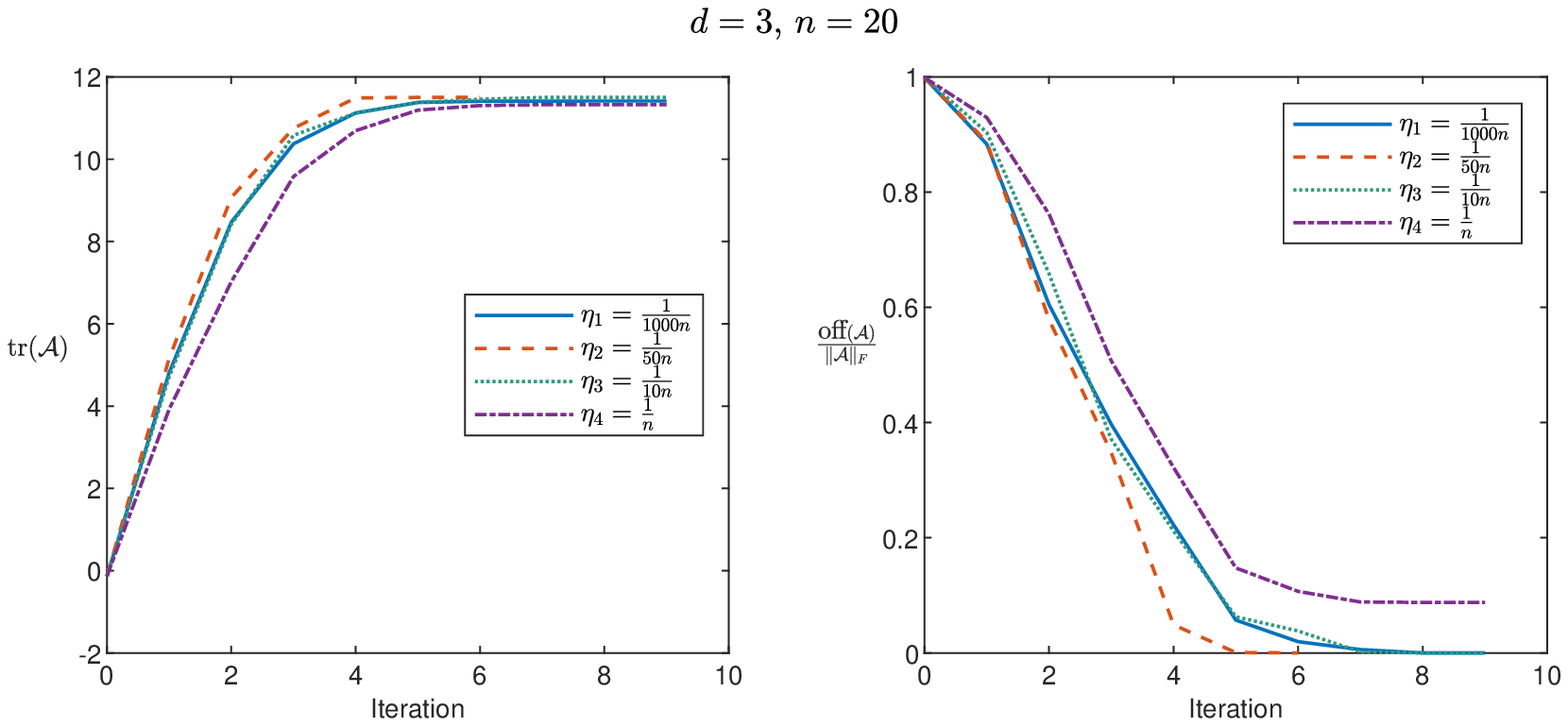} \\
    \includegraphics[width=\textwidth]{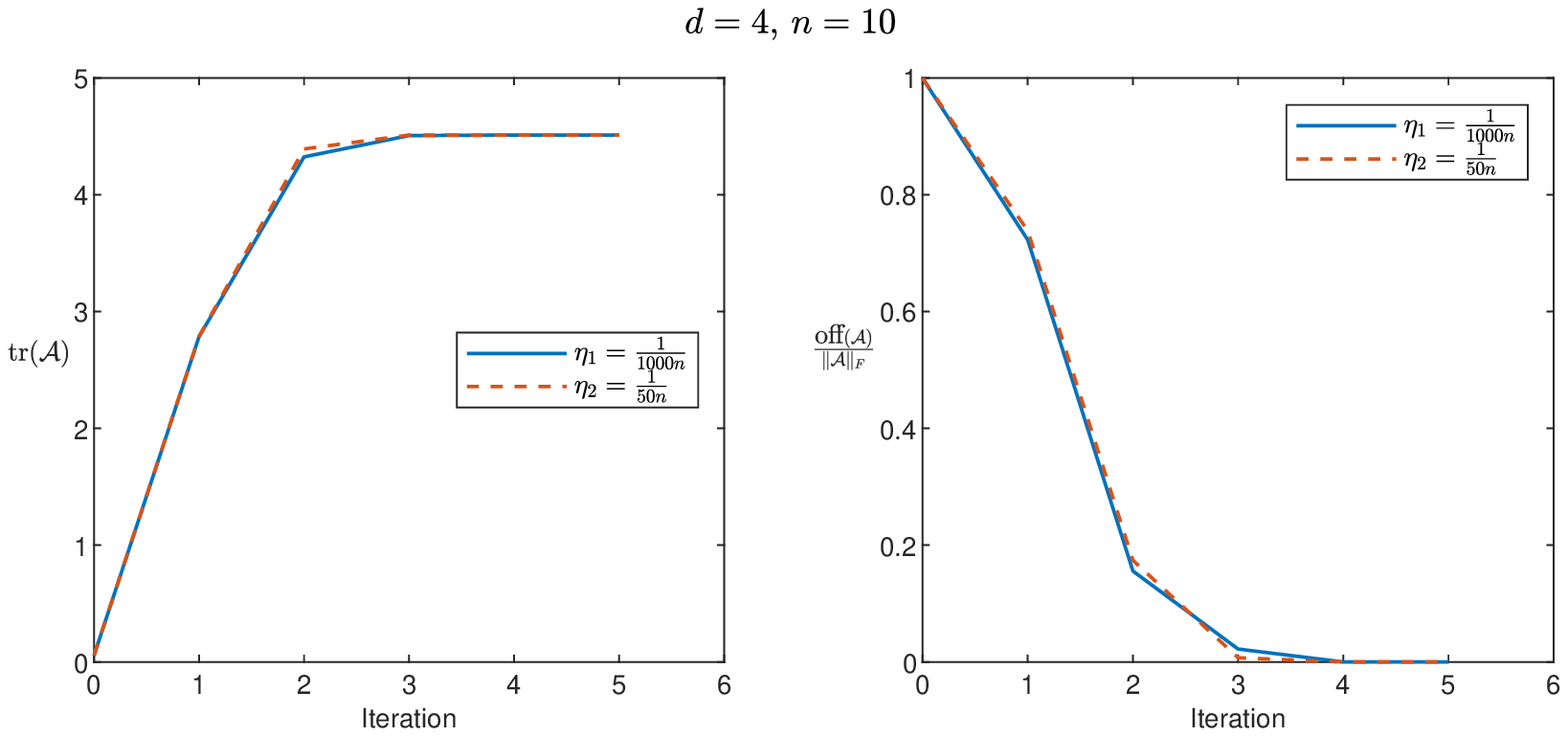}
    \caption{Convergence of Algorithm~\ref{agm:jacobi} for different values of $\eta$ on tensors of order $3$ and $4$ that are diagonalizable using orthogonal transformations.}\label{fig:diag}
\end{figure}

Figure~\ref{fig:diag} shows the convergence of the trace and the relative off-norm in the Algorithm~\ref{agm:jacobi} for diagonalizable $20\times20\times20$ and $10\times10\times10\times10$ tensors, for different values of $\eta$ from~\eqref{condition}. One can observe that for a high value of $\eta$, $\eta=\frac{1}{n}$, the trace converges to a lower value than for the smaller $\eta$. Moreover, in these examples for $\eta=\frac{1}{n}$ the relative off-norm converges to a number greater than $0$, while for smaller $\eta$ it converges to $0$. This means that for $\eta=\frac{1}{n}$ the algorithm converges to a different stationary point than for smaller $\eta$, the one that is not a diagonal tensor. Therefore, from our observations we recommend the use of smaller $\eta$.

\begin{figure}[h!]
    \centering
    \includegraphics[width=\textwidth]{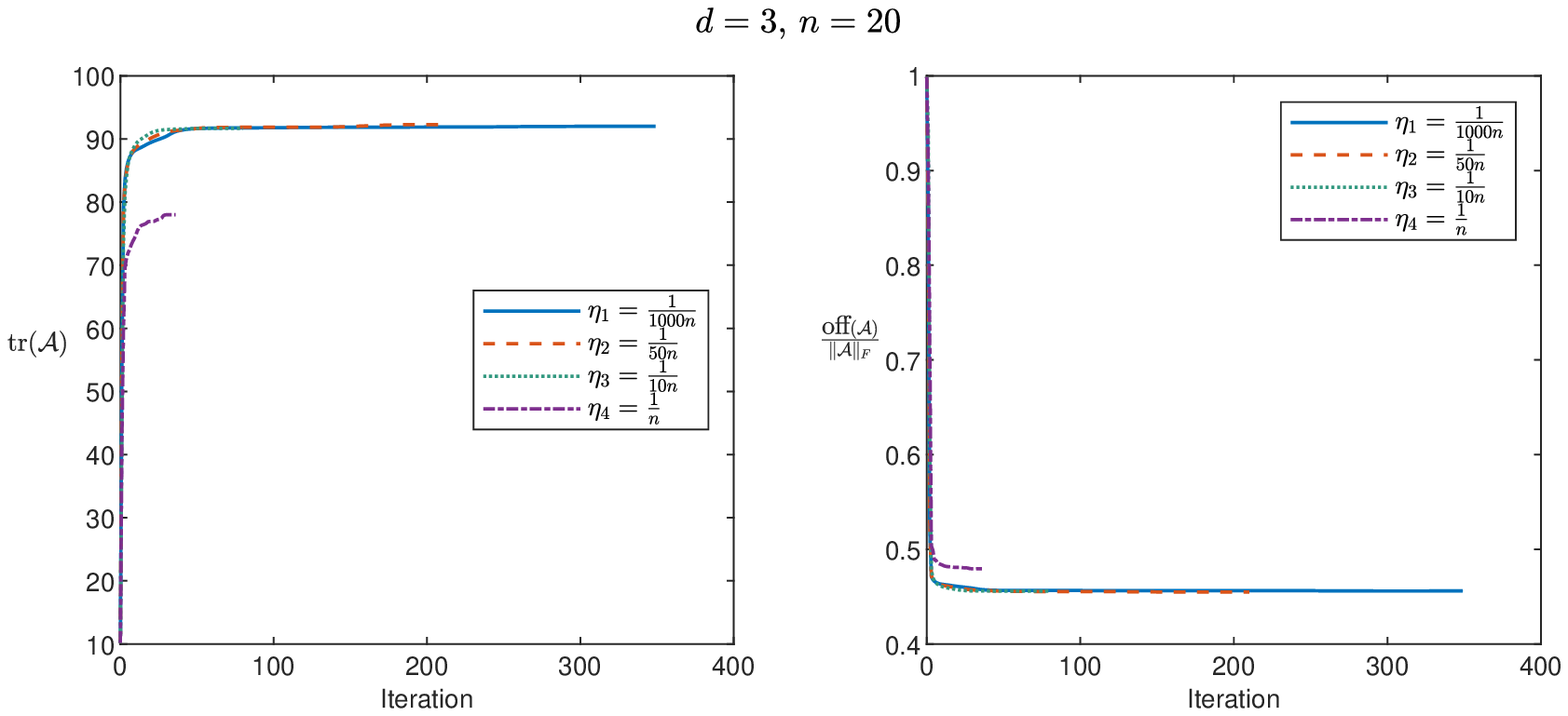} \\
    \includegraphics[width=\textwidth]{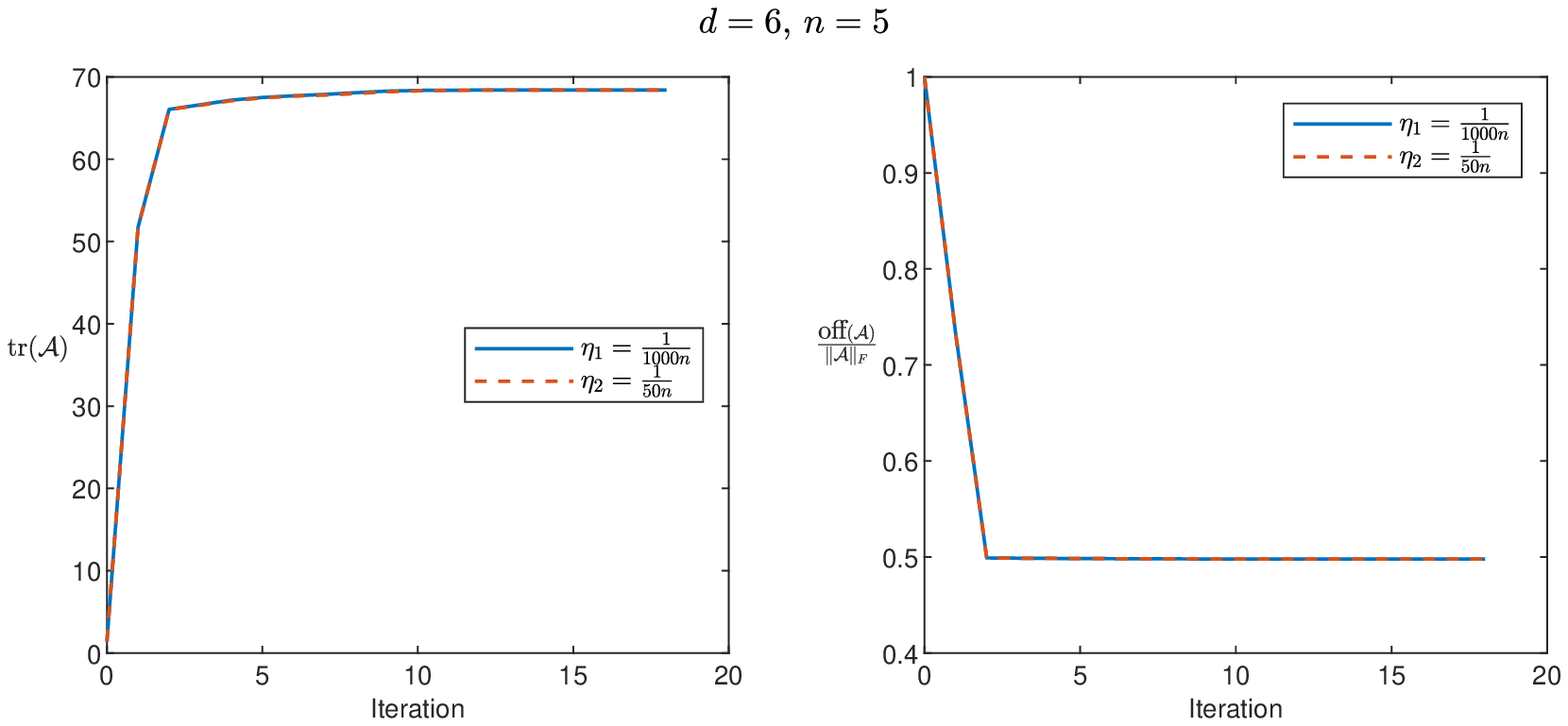}
    \caption{Convergence of Algorithm~\ref{agm:jacobi} for different values of $\eta$ on random tensors of order $3$ and $6$.}\label{fig:nondiag}
\end{figure}

We repeat the same experiment as the one described above, but this time on non-diagonalizable $20\times20\times20$ and $5\times5\times5\times5\times5\times5$ tensors. Here one cannot expect the relative off-norm to be equal to zero. The results are shown in Figure~\ref{fig:nondiag}. Same as in Figure~\ref{fig:diag}, for $\eta=\frac{1}{n}$ we get the convergence to a different, less desirable, stationary point of the function $f$.

\begin{figure}[h!]
    \centering
    \includegraphics[width=0.45\textwidth]{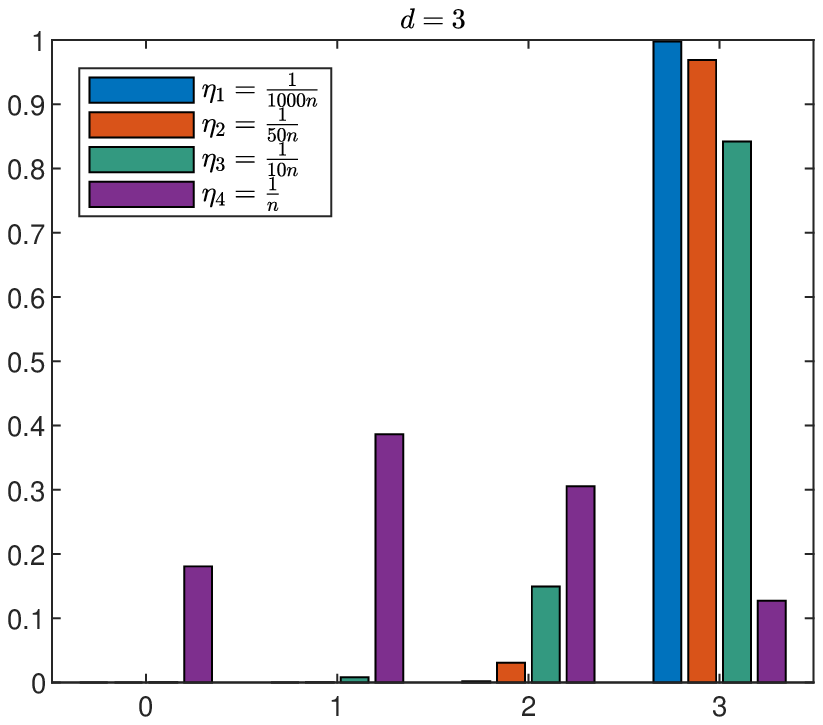}
    \includegraphics[width=0.45\textwidth]{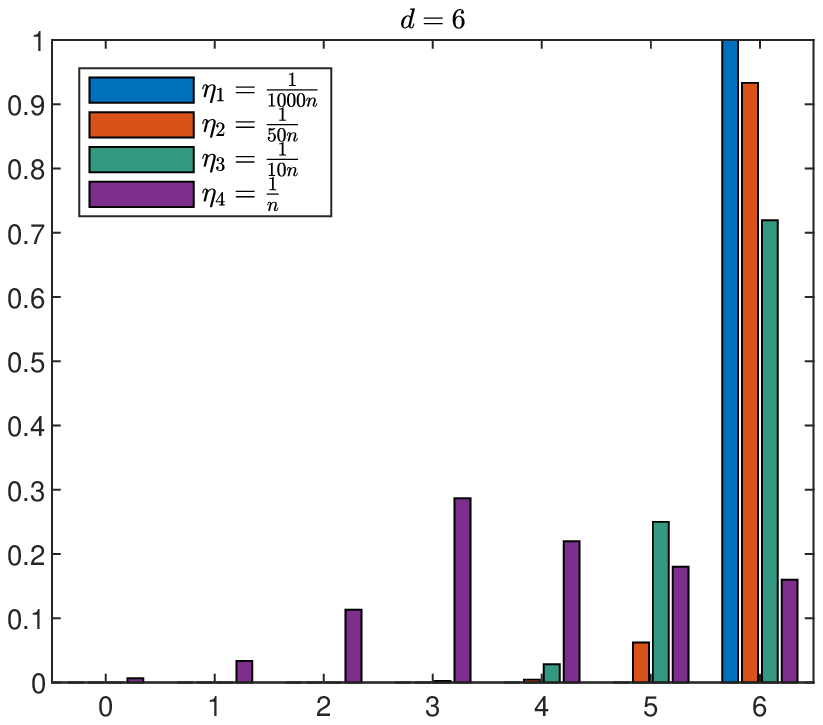}
    \caption{Portion of the number of microiterations within one iteration for different values of $\eta$ on random tensors of order $3$ and $6$.}\label{fig:microit}
\end{figure}

Barplots in Figure~\ref{fig:microit} show how $\eta$ affects the number of microiterations in each iteration of the Algorithm~\eqref{agm:jacobi}. The test is done on non-diagonalizable $20\times20\times20$ and $5\times5\times5\times5\times5\times5$ tensors. If $\eta$ is bigger, the condition~\eqref{condition} is more restrictive and more microiterations are skipped. For example, for $d=3$ and $\eta=\frac{1}{n}$, $38.6\%$ of iterations contain only one microiteration, and only $12.7\%$ contain the maximal number of microiterations. On the other hand, for $\eta=\frac{1}{1000n}$, $99.8\%$ of iterations consists of three microiterations.

\begin{figure}[h!]
    \centering
    \includegraphics[width=0.45\textwidth]{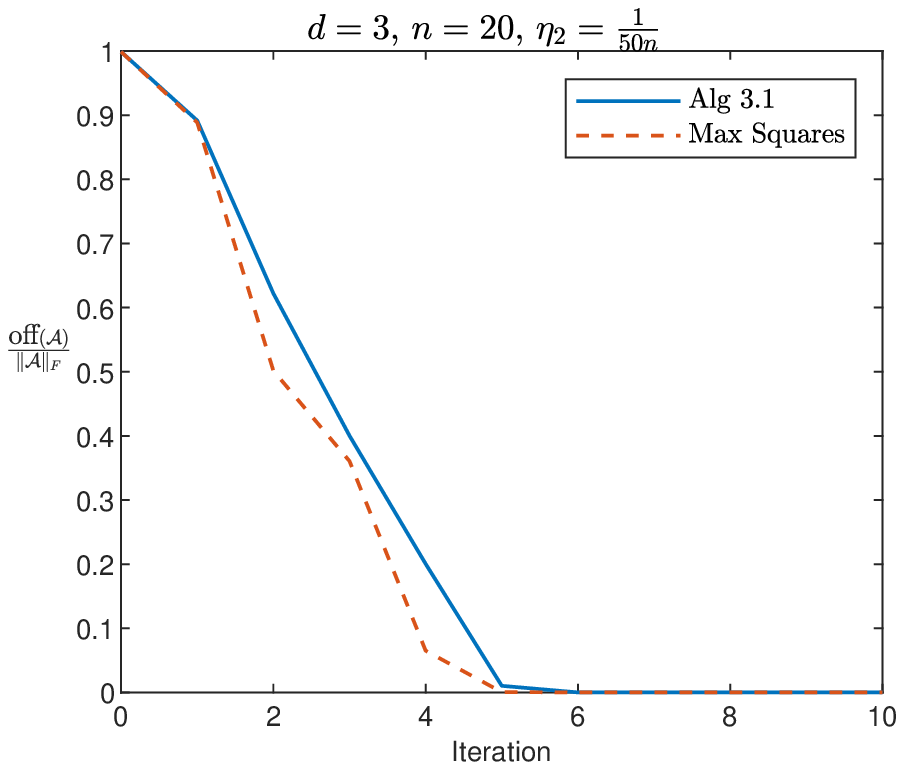}
    \includegraphics[width=0.45\textwidth]{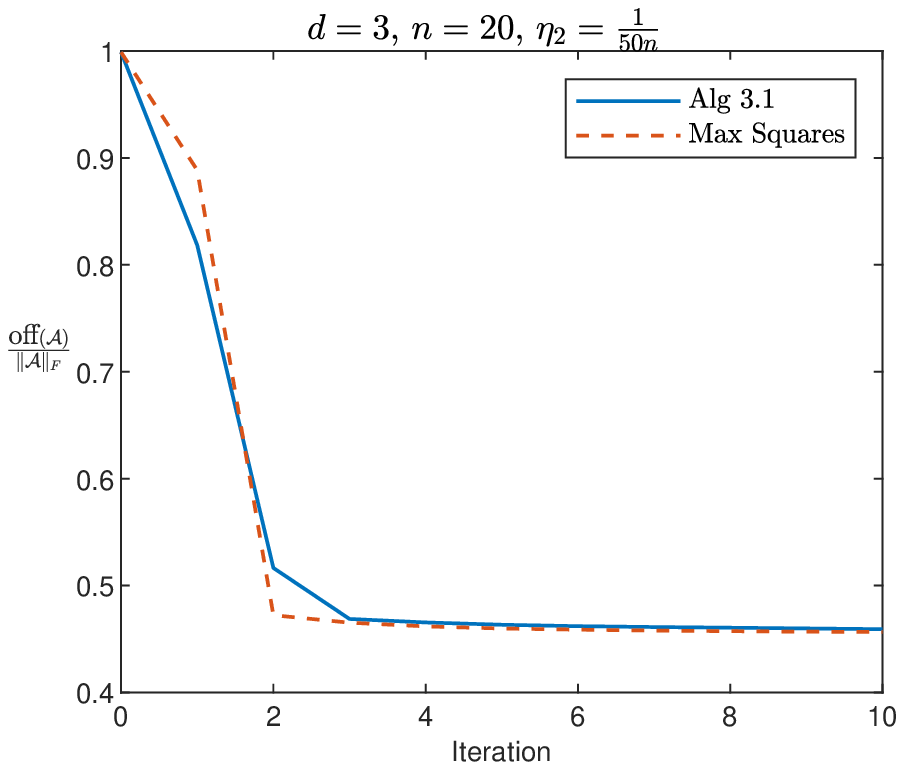}
    \caption{Trace maximization compared to the maximization of the squares of the diagonal elements for two random tensors during the first $10$ iterations.}\label{fig:tracevssquares}
\end{figure}

In Figure~\ref{fig:tracevssquares} we consider the trace maximization algorithm opposed to the Jacobi-type algorithm that maximizes the squares of the diagonal elements. We observe the performance of the Algorithm~\ref{agm:jacobi} and the algorithm from~\cite{Be21} on two random $20\times20\times20$ tensors, one orthogonally diagonalizable, and one non-diagonalizable. We can see that the results of both algorithms are comparable.

\begin{figure}[h!]
    \centering
    \includegraphics[width=\textwidth]{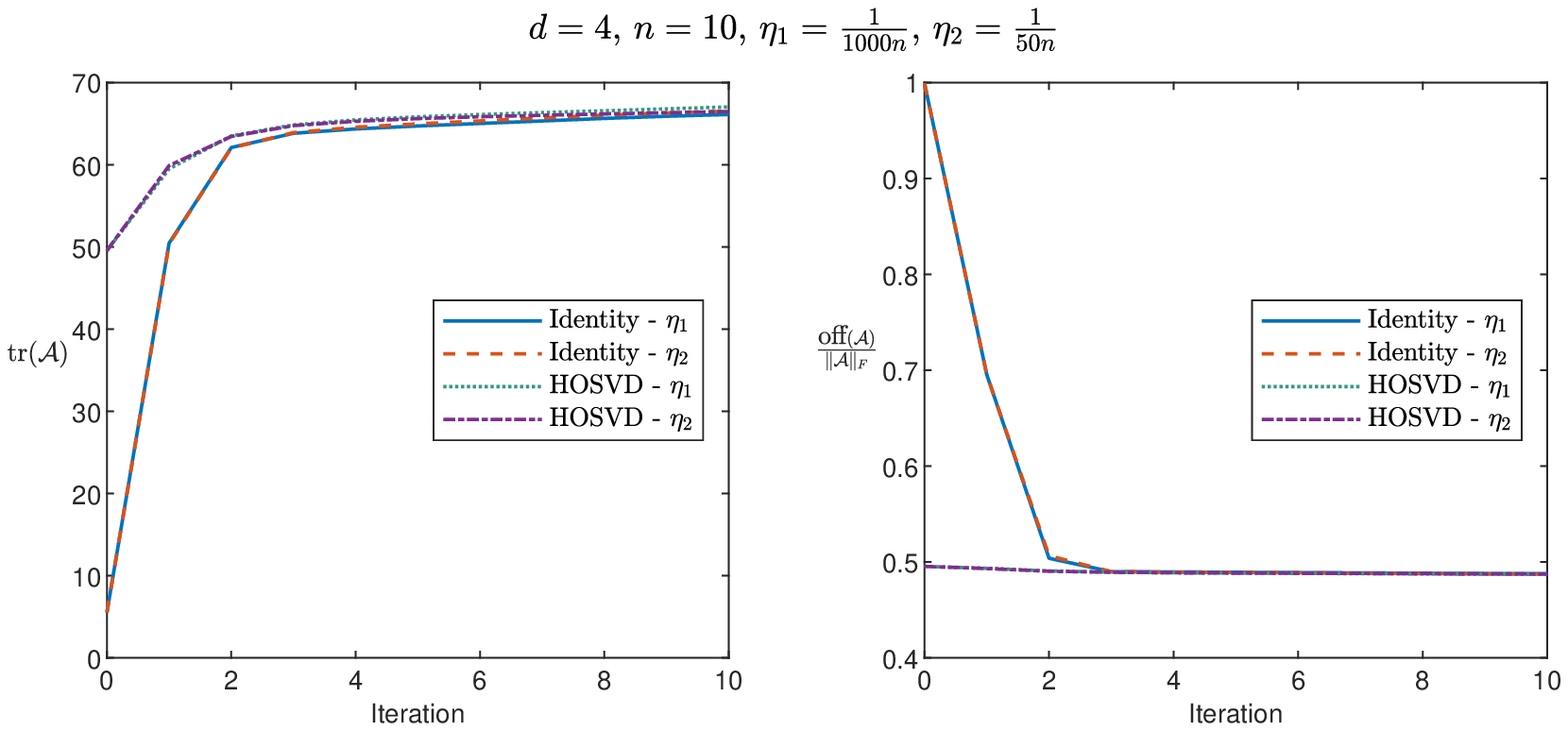}
    \caption{Comparison of different initialization strategies for Algorithm~\ref{agm:jacobi} for a random tensor of order $4$ during the first $10$ iterations.}\label{fig:hosvd}
\end{figure}

In Section~\ref{sec:algorithm} we discussed different initialization strategies for the Algorithm~\ref{agm:jacobi}. In Figure~\ref{fig:hosvd} we compare the identity initialization given by~\eqref{identinit} and the HOSVD initialization given by~\eqref{hosvdinit}. The results are shown for a non-diagnalizable $10\times10\times10\times10$ tensor during the first $10$ iterations. When a tensor is preconditioned using HOSVD it becomes closer to a diagonal one. Thus, the starting trace value is higher for the HOSVD initialization than for the identity initialization. Also, the starting relative off-norm is much closer to the limit value, significantly under $1$. Regardless, after the first few iterations the algorithm using identity initialization catches up. Both initializations give equally good approximations and converge to the same value.

\begin{figure}[h!]
    \centering
    \includegraphics[width=\textwidth]{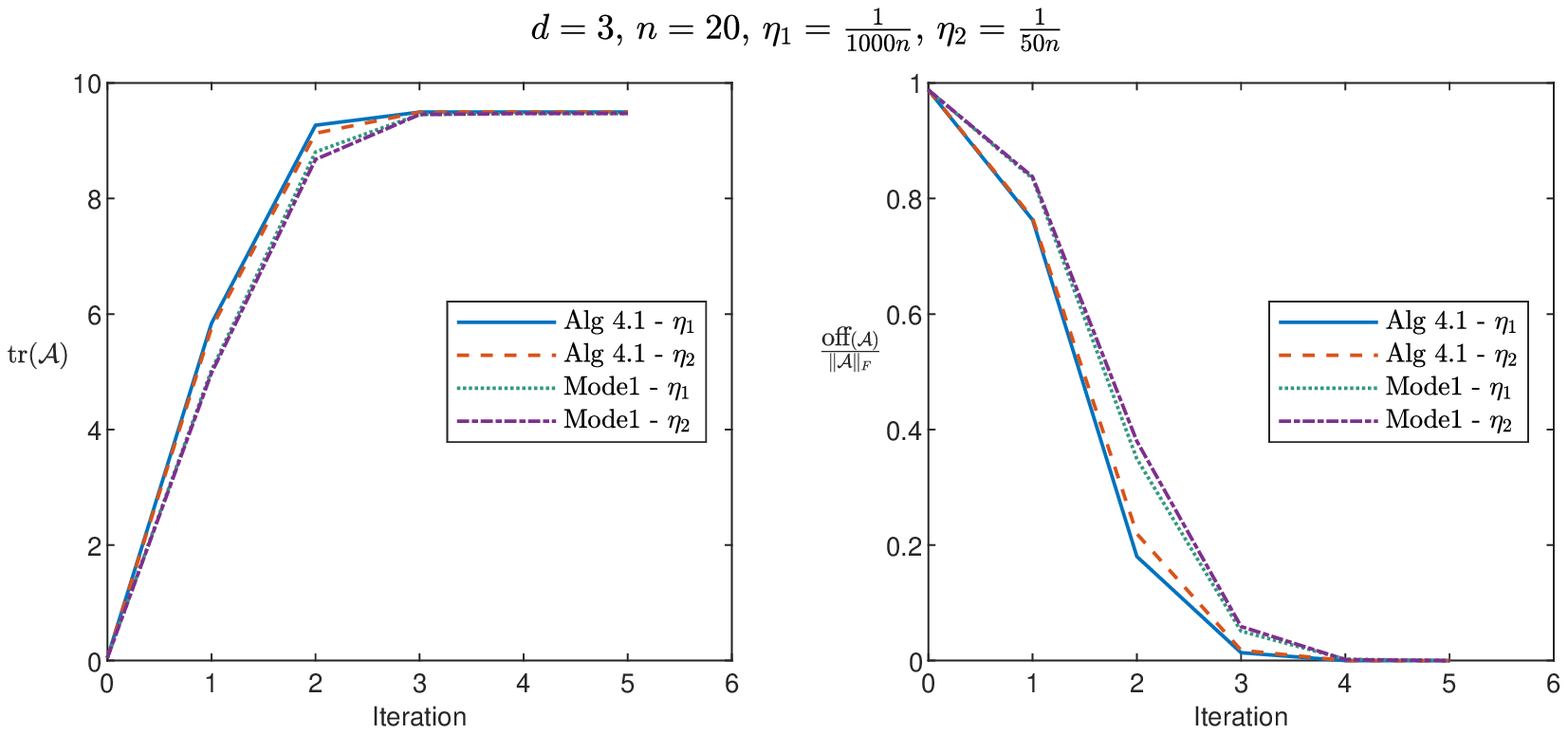}
    \caption{Convergence of Algorithm~\ref{agm:symjacobi} and its ``Mode1'' modification for different values of $\eta$ on a tensor of order $3$ that is diagonalizable using orthogonal transformations.}\label{fig:symdiag}
\end{figure}

Lastly, we observe the symmetric case. We present the convergence results of Algorithm~\ref{agm:symjacobi} for a diagonalizable $20\times20\times20$ tensor in Figure~\ref{fig:symdiag}.

Additionally to the Algorithm~\ref{agm:symjacobi} we observe its modification, a hybrid approach between the Algorithms~\ref{agm:jacobi} and~\ref{agm:symjacobi}. In the Algorithm~\ref{agm:symjacobi} the rotation matrix is chosen by optimizing the angle with respect to all modes, which leads to a cubic equation~\eqref{symtanphi} for $d=3$. As this can be computationally challenging, we investigated another approach. From the relation~\eqref{tanphi} used in the Algorithm~\ref{agm:jacobi} we compute the rotation angle that is optimal in one mode $l$, e.g. $l=1$. Because of the symmetry, it does not matter which mode we consider. Then we apply this same rotation in all modes, as it is done in the Algorithm~\ref{agm:symjacobi}, to preserve the symmetry. This modification is denoted by \emph{Mode1} in Figure~\ref{fig:symdiag}.
The convergence results from Section~\ref{sec:cvg} hold for Algorithm~\ref{agm:symjacobi}, but not for the modification \emph{Mode1}. In practice we observed that \emph{Mode1} converges to the same solution as Algorithm~\ref{agm:symjacobi}, only slower. Analogue results were obtained for the non-diagonalizable tensors, too.

\begin{figure}[h!]
    \centering
    \includegraphics[width=\textwidth]{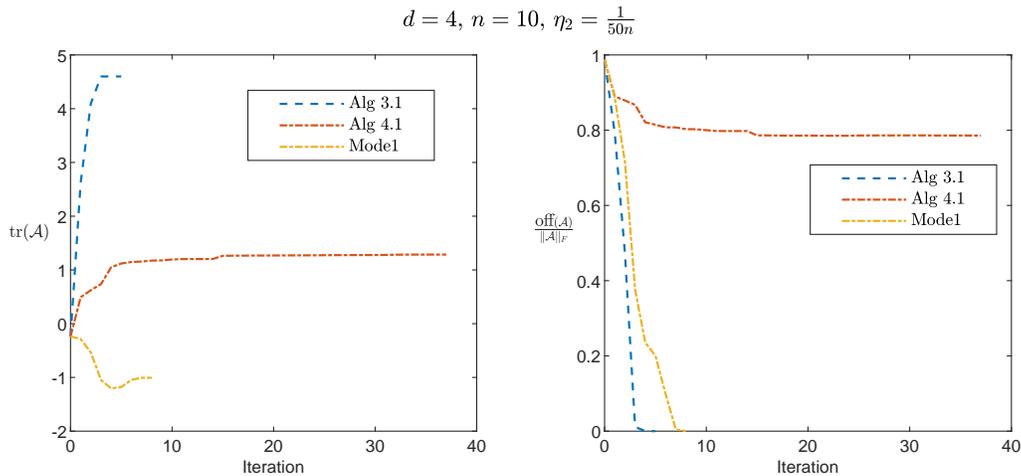}
    \caption{Convergence of Algorithm~\ref{agm:jacobi},  Algorithm~\ref{agm:symjacobi}, and its ``Mode1'' modification on a symmetric 4th-order tensor that is diagonalizable using orthogonal transformations.}\label{fig:symdiag4d}
\end{figure}

Still, the symmetry-preserving Algorithm~\ref{agm:symjacobi} has some limitations when dealing with tensors of even order, specifically tensors with both positive and negative elements on the diagonal approximation~\cite{SLCID12}. Although the convergence theorem is still valid, the acquired stationary point of the objective function is not its global maximum. This behaviour can be seen in Figure~\ref{fig:symdiag4d} where we apply the Algorithm~\ref{agm:jacobi}, the Algorithm~\ref{agm:symjacobi}, as well as its Mode1 modification on a random orthogonally diagonalizable symmetric 4th-order tensor. A symmetric 4th-order diagonalizable tensor is constructed by taking a diagonal tensor with diagonal elements drawn uniformly from $[-1,1]$, that is multiplied in each mode by a same random orthogonal matrix. In each iteration of the Algorithm~\ref{agm:symjacobi}, the rotation matrix is chosen by solving the 4th-order equation (analogue to equation~\eqref{symtanphi} for $d=3$),
\begin{align*}
& (a_{1222}-a_{1112})t^4-(a_{1111}+a_{2222}-6a_{1122})t^3 -6(a_{1222}-a_{1112})t^2+ \nonumber \\
& \qquad\qquad\qquad\qquad\quad+(a_{1111}+a_{2222}-62a_{1122})t+(a_{1222}-a_{1112})=0.
\end{align*}
In numerical experiments we observed some interesting things. Algorithm~\ref{agm:jacobi} yielded matrices $U_l$, $1\leq l\leq 4$, equal up to sign, specifically
$$-U_1=U_2=U_3=U_4.$$
Moreover, obtained diagonal elements are good approximations of the absolute values of eigenvalues of the starting tensor, except for one value which is of negative sign.
Next, Mode1 modification does not find the maximal trace. However, the approximation is as good as for the Algorithm~\ref{agm:jacobi} (in terms of the off-norm). Additionally, obtained diagonal elements are good approximations of the eigenvalues of the starting tensor.

\begin{thebibliography}{10}
\bibitem{Be21}
E. Begovi\'{c}~Kova\v{c}:
\emph{Convergence of a Jacobi-type method for the approximate orthogonal tensor diagonalization}.
preprint, arXiv:2109.03722 [math.NA]
%
\bibitem{Cich15}
A. Cichocki, D. Mandic, C. Caiafa, A.-H. Phan, G. Zhou, Q. Zhao, L. De Lathauwer:
\emph{Multiway component analysis: tensor decompositions for signal processing applications}.
IEEE Sig. Process. Mag. 32(2) (2015) 145--163.
%
\bibitem{ComonSP}
P. Comon:
\emph{Tensor diagonalization, a useful tool in signal processing}.
IFAC Proceedings Volumes 27(8) (1994) 77--82.
%
\bibitem{ComSor07}
P. Comon, M. Sorensen:
\emph{Tensor diagonalization by orthogonal transforms}.
Research Report ISRN I3S/RR-2007-06-FR (2007).
%
\bibitem{DeLHOSVD}
L. De Lathauwer, B. De Moor, J. Vandewalle:
\emph{A Multilinear Singular Value Decomposition}.
SIAM J. Matrix Anal. Appl. 21(4) (2000) 1253--1278.
%
\bibitem{DeLDeMV2001}
L. De Lathauwer, B. De Moor, J. Vandewalle:
\emph{Independent component analysis and (simultaneous) third-order tensor diagonalization}.
IEEE Transactions on Signal Processing 49 (2001) 2262--2271.
%
\bibitem{Ishteva13}
M. Ishteva, P.-A. Absil, P. Van Dooren:
\emph{Jacobi algorithm for the best low multilinear rank approximation of symmetric tensors}.
SIAM J. Matrix Anal. Appl. 34(2) (2013) 651--672.
%
\bibitem{KB09}
T.~G. Kolda, B.~W. Bader:
\emph{Tensor decompositions and applications}.
SIAM Rev. 51(3) (2009) 455--500.
%
\bibitem{LUC19}
J. Li, K. Usevich, P. Comon:
\emph{On approximate diagonalization of third order symmetric tensors by orthogonal transformations}.
Linear Algebra Appl. 576(1) (2019) 324--351.
%
\bibitem{LUC18}
J. Li, K. Usevich, P. Comon:
\emph{Globally Convergent Jacobi-Type Algorithms for Simultaneous Orthogonal Symmetric Tensor Diagonalization}.
SIAM J. Matrix Anal. Appl. 39(1) (2018) 1--22.
%
\bibitem{MMVL08}
C. D. Moravitz Martin, C. F. Van Loan:
\emph{A Jacobi-type method for computing orthogonal tensor decompositions}.
SIAM J. Matrix Anal. Appl. 30(3) (2008) 1219--1232.
%
\bibitem{PTC17}
A. H. Phan, P. Tichavsk\'y, A. Cichocki:
\emph{Blind Source Separation of Single Channel Mixture Using Tensorization and Tensor Diagonalization.}
In: P. Tichavsk\'y, M. Babaie-Zadeh, O. Michel, N. Thirion-Moreau (eds) Latent Variable Analysis and Signal Separation. LVA/ICA 2017. Lecture Notes in Computer Science, vol 10169. Springer, 2017.
%
\bibitem{SLCID12}
M. Sorensen, L. De Lathauwer, P. Comon P, S. Icart, L. Deneire:
\emph{Canonical Polyadic decomposition with a Columnwise Orthonormal Factor Matrix}.
SIAM J. Matrix Anal. Appl. 33(4) (2012) 1190--1213.
%
\bibitem{TPC17}
P. Tichavsk\'y, A. H. Phan, A. Cichocki:
\emph{Non-orthogonal tensor diagonalization}.
Signal Process. 138 (2017) 313--320.
%
\bibitem{LUC20}
K. Usevich, J. Li, P. Comon:
\emph{Approximate matrix and tensor diagonalization by unitary transformations: convergence of Jacobi-type algorithms}.
SIAM J. Optim. 30(4) (2020) 2998--3028.
\end{thebibliography}
\end{document}